\documentclass[11pt]{amsart}
\usepackage{amssymb}
\newcommand{\BP}{{\mathbb{P}}}
\newcommand{\BA}{{\mathbb{A}}}
\newcommand{\gb}{\beta}
\newcommand{\ga}{\alpha}
\newcommand{\gc}{\gamma}
\newcommand{\Pic}{\mathrm{Pic}}
\newcommand{\ra}{\rightarrow}
\newcommand{\ol}[1]{\overline{#1}}
\swapnumbers
\theoremstyle{plain}
\newtheorem{lma}{Lemma}
\newtheorem{thm}[lma]{Theorem}

\theoremstyle{definition}
\newtheorem{dfn}[lma]{Definition}
\newtheorem{rmr}[lma]{Remark}
\newtheorem{ntt}[lma]{Notation}

\newtheorem{dsc}[lma]{}

\begin{document}

\title{The chord construction}
\author{David Lehavi}
\address{Mathematics Institute\\
The Hebrew University of Jerusalem\\
Givat Ram, Jerusalem  91904,  Israel}
\email{dlehavi@math.huji.ac.il}
\author{Ron Livn\'e}
\address{Mathematics Institute\\
The Hebrew University of Jerusalem\\
Givat Ram, Jerusalem  91904,  Israel}
\email{rlivne@math.huji.ac.il}
\date{\today}
\keywords{Plane cubics, Chords}
\subjclass{14N15,14H52}
\begin{abstract}
Let $F$ be a smooth plane curve of degree 3. Let $\gb$ be an element in
$\Pic(F)[2]-\{0\}$. Let us define
\[
  F':=\{\ol{p(p+\gb)}|p\in F\}\subset(\BP^2)^*.
\]
In this note we show that $F'$ is a smooth embedding of $F/\gb$.
Moreover, let $\gb'$ be the generator of
$\Pic(F)/\gb$, and let $p\in F$ be a flex, then $\ol{p(p+\gb)}+\gb'$ is a
flex on $F'$. We present two proofs.
\end{abstract}
\maketitle

\begin{ntt}
Let $F$  be a smooth curve of genus 1, $\gb\in\Pic(F)[2]$. Let us denote
$F':=F/\gb$. Let $\pi$ be the quotient map $F\ra F'$ and let $\gb'$ be the
image of the non-zero element in $\Pic(F)[2]/\gb$ in $\Pic(F')$.
Let $H$ be an element in $\Pic^3(F)$, and $i:F\ra\BP^2$ the map associated to
$H$.
This notation will hold in the next definition and couple of lemmas.
\end{ntt}
\begin{dfn}[The chord construction (c.c.)]\label{Dcc}
Let us define
\[
  \begin{aligned}
  i':F/\gb  &\ra{\BP^2}^* \\
   \{p,p+\gb\}&\mapsto\ol{i(p)i(p+\gb)}
  \end{aligned}
\]
\end{dfn}
\begin{thm}\label{Lccembed}
The $c.c$ is an embedding.
\end{thm}
\begin{proof}
Let us identify $F$ with its image $i(F)$.
\item {\em The c.c. is 1-1:} Let $H$ be the line corresponding to a point
$\{p,\gb+p\}\in F/\gb$. Then  $F\cap H$ consists of three points:
\[
  \{p,p+\gb,H-(2p+\gb)\}.
\]
If $(H-(2p+\gb))-\gb\neq p,p+\gb$, then $i'(H)$ consist of the unique point
$\{p,p+\gb\}$, so it is 1-1.
If $(H-(2p+\gb))-\gb=p$ then $H-(2p+\gb)=p+\gb$, so $i'$ is 1-1 here too.
The remaining case is similar.
(in the last two cases $p$ (respectively $p+\gb$) is a flex).
For $p\in\BP^2$, let $p^*$ be the line it defines in ${\BP^2}^*$
\item {\em If $F$ does not have a CM then for a generic $p\in F$, $p^*$ is transversal to $F'$:}
Suppose it is not.
By upper-semicontinuity, it is never transversal, so $F$ is the dual curve
of $i'(F')$. Since $i'(F')$ is not a line (it is an immersion of $F'$),
this induces a degree 1 map $i(F')\ra F$, and therefor a degree 1 map
$F'\ra F$. This means that $F$ has a CM.
\item {\em The c.c is of degree 3:}
Assume first that $F$ does not have a CM.
Choose a generic $p\in F$. The degree of $i'$ is the number of lines
through $p$ which belong to $i'(F')$.
There are exactly two types of such lines:
\begin{itemize}
  \item  The line is $\ol{p(p+\gb)}$.
  \item  The line is the one passing through a pair of the four solutions of
     $x+(x+\gb)=H-p$.
\end{itemize}
Now suppose $F$ has a CM. Then $F$ and the c.c. are a flat limit of curves
without a CM, so the degree of $i'$ is constant. (see \cite{Ha} preposition 
III.9.8)
\end{proof}
\begin{thm}\label{Lccflex}
If $i(p)$ is a flex, then so is $i'(\pi(p)+\gb')$.
\end{thm}
\begin{proof}
We are going to drop the $i,i'$ notation, and identify both curves with their
images.
\item {\em Involutions on $F\subset\BP^2$:} 
Let $p$ be a flex on $F$. Let $\ga$ be the unique involution of $\BP^2$ that
\begin{itemize}
  \item sends $F$ to itself.
  \item fixes $p$.
\end{itemize}
Note that this involution has a fixed line $L$, where $L\cap F$ consists of the
translates of $p$ by points in $\Pic(F)[2]-\{0\}$. All the involutions of
$\BP^2$ that send $F$ to itself are of this form (they have to send flexes
to flexes, and there are 9 flexes, so one of them is fixed).
The same holds for $F'$.
\item{\em The dual involution:} Given $\ga$, there is a dual involution $\ga'$
on ${\BP^2}^*$. By the definition of the c.c, $\ga'$ sends $F'$ to itself.
The fixed line and point of $\ga'$ are $p^*,L^*$ respectively.
There is exactly one pair of points in $F$ with difference $\gb$ on $L$.
This is the pair $\{p+\gc,p+\gc+\gb\}$ where $\gc$ is some element
in $\Pic(F)[2]-\{0,\gb\}$.
\end{proof}
\begin{dsc}
We develop a formula for the c.c. . Theorem \ref{Tccformula}
(which states the formula) also gives an alternative proof for theorems
\ref{Lccembed},\ref{Lccflex}.
\end{dsc}
\begin{ntt}
Let us denote homogeneous coordinates on $\BP^2$ by $[X;Y;Z]$. Let
us denote the dual coordinates on ${\BP^2}^*$ by $[U;V;W]$.
Let us also define
\[
  x:=X/Z, y:=Y/Z,
\]
and denote the $\BA^2$ coordinates by $(x,y)$
Let $F$ as the closure of the locus of $y^2=f$ where
\[
  f(x)=x^3+ax^2+bx.
\]
Let us identify $F$ with $\Pic(F)$ by choosing $O:=[0;1;0]$ be the zero.
Let $\gb$ be the point $(0,0)=[0;0;1]$.
\end{ntt}
\begin{lma}
Let $p=(x(p),y(p))$ be a point on $F$.
Then
\[
  p+\gb=(bx(p),-by(p)/x(p)^2)
\]
\end{lma}
\begin{proof}
There are exactly three points in
\[
  F\cap\ol{\gb p},
\]
these are $p,\gb,-p-\gb$.
Since $\ol{\gb p}$ is given by the nulls of $y/x=\frac{y(p)}{x(p)}$,
They are the solutions of
\[
  (x\frac{y(p)}{x(p)})^2=x^3+ax^2+bx.
\]
One of the solutions of this equation is $0$.
the other two are the roots of
\[
  0=x^2+(a-(\frac{y(p)}{x(p)})^2)x+b.
\]
Since one of the roots is $x(p)$, we see that the other is
\[
  b/x(p)=x(p+\gb).
\]
We also deduce that
\[
  y(p+\gb)=-y(-p-\gb)=-\frac{y(p)}{x(p)}b/x(p)=-by(p)/x(p)^2.
\]
\end{proof}
\begin{lma}
The lines of the c.c. are given in ${\BP^2}^*$ by
\[
  [y(p)(x(p)^2+b);bx(p)-x(p)^3;-2bx(p)y(p)],
\]
where $p$ are the points on $F$.
\end{lma}
\begin{proof}
Suppose the line $\ol{p(\gb+p)}$ is given by
\begin{equation}\label{Esxcy}
  sx+c=y.
\end{equation}
In this case we have:
\[
  \begin{aligned}
  s&=\frac{y(p)-y(p+\gb)}{x(p)-x(p+\gb)}
    =\frac{y(p)+by(p)/x(p)^2}{x(p)-b/x(p)}=\frac{y(p)(x(p)^2+b)}{x(p)^3-bx(p)}
    \\
  C&=y(p)-sx(P)
  \end{aligned}
\]
Simplifying
\[
  [s;-1;c],
\]
We get the result.
\end{proof}
\begin{lma}
Assuming the c.c. is of degree 3, it satisfies the equation
\[
e(U-\mu^{-1}W)V^2 = W^3-\frac{2\mu b}{\mu^2-b} (U-\mu^{-1}W)W^2 - 
                    \frac{\mu^2 b}{\mu^2-b} (U-\mu^{-1}W)^2 W.
\]
for some $e$, where $\mu=2b/a$.
\end{lma}
\begin{proof}
Let us consider the involution $(x,y)\mapsto(x,-y)$ on $\BP^2$.
It has a fixed line $l$: the nulls of $y=0$, and a unique fixed point out
of $l$: $O$.
If we consider the dual involution, we see that it has a fixed line at the
dual of $O$, to be denoted $O^*$. It also has a unique fixed point out of $O^*$
. This is the point dual to $l$, to be denoted $l^*$.
This means that $l^*$ is a flex of the c.c. and its translations by elements
of order 2, all sit on $O^*$.
As in the proof of lemma \ref{Lccembed}, there are two types of such lines:
\begin{itemize}
  \item The line is $\ol{O\gb}$.
  \item The line is the one passing through a pair of the four solutions of
     $2x=\gb$.
\end{itemize}
In the first type, this is the line $x=0$.
In the second type, we are looking for the solutions of
\[
  \begin{aligned}
   x&\neq 0,y^2=f(x),x\frac{dy}{dx}=y &\Longleftrightarrow \\
   x&\neq 0,x\frac{df}{dx}=f &\Longleftrightarrow \\
   3&x^2+2ax+b=2x^2+2ax+2b &\Longleftrightarrow \\
   x&^2=b.
  \end{aligned}
\]
We have to identify the tangent to the flex.
A tangent to the flex (on a curve of degree 3) is a line that intersect the
curve only in the flex.
In dual coordinates:
It is the unique point $q=(\mu,0)$ on the line $y=0$, such that $y=0$,
is the ONLY c.c line on which $q$ sits.
The $x$ coordinates of the intersection between a chord and $y=0$, can be
obtained from equation \ref{Esxcy}, letting $y=0$. i.e. they are:
\[
  -c/s=\frac{y(p)-sx(p)}{s}=y(p)/s-x(p)
  =\frac{x(p)^3-bx(p)}{x(p)^2+b}-x(p)
  =\frac{-2bx(p)}{x(p)^2+b}.
\]
We see that each such point is given by only two $x(P)$.
In $q$, these $x(p)$ are given by the roots of
\[
  x^2+ax+b.
\]
Whence,
\[
  \mu=2b/a.
\]
Let us apply the projective transformation
\[
U'=U-1/\mu W,V'=V,W'=W,
\]
and define
\[
 w'=W'/U',v'=V'/U'.
\]
In these coordinates we have
\[
\begin{aligned}
e {v'}^2 & =
 w'(w'-\mu\sqrt{b}/(\mu-\sqrt{b})) (w'+\mu\sqrt{b}/(\mu+\sqrt{b})) \\
 & =w'({w'}^2-\frac{2\mu b}{\mu^2-b}w'-\frac{\mu^2 b}{\mu^2-b})
\end{aligned}
\]
for some e. Going back to our original coordinates, we get
\[
e (U-\mu^{-1}W)V^2 = W^3-\frac{2\mu b}{\mu^2-b} (U-\mu^{-1}W)W^2 - 
                    \frac{\mu^2 b}{\mu^2-b} (U-\mu^{-1}W)^2 W.
\]
\end{proof}
\begin{thm}\label{Tccformula}
The c.c is of degree 3, and $e=\frac{-8b^3}{a^2-4b}$
\end{thm}
\begin{proof}
This is a trivial check, best done by a computer.
One has to define:
\[
  U:=y(x^2+b),V:=bx-x^3,W:=-2bxy,\mu:=2b/a,
\]
and verify that
\[
  \frac{W^3-\frac{2\mu b}{\mu^2-b} (U-\mu^{-1}W)W^2 -
  \frac{\mu^2 b}{\mu^2-b} (U-\mu^{-1}W)^2 W
  }
  {(U-\mu^{-1}W)V^2}
  =-\frac{8b^3 y^2}{(a^2-4b)f(x)}
\]
\end{proof}
\begin{rmr}
If $F$ is as before, and $\gb\in\Pic^0(F)-\Pic^0(F)[2]$, the chord construction
gives a map of degree 1. The image of this map is of degree 6.
\end{rmr}


\begin{thebibliography}{99}
\bibitem[Ha]{Ha}  R. Hartshorne {\em Algebraic Geometry}
              Springer Verlag, GTM 52. (1977)
\end{thebibliography}
\end{document}